\documentclass[11pt]{amsart}

\theoremstyle{plain}
\newtheorem{thm}{Theorem}[section]
\newtheorem{theorem}[thm]{Theorem}

\newtheorem{proposition}[thm]{Proposition}
\theoremstyle{definition}
\newtheorem{remark}[thm]{Remark}

\newtheorem{definition}[thm]{Definition}

\newtheorem{question}[thm]{Question}

\numberwithin{equation}{section}

\newcommand{\wq}{\widetilde{q}}
\newcommand{\wt}{\widetilde{\theta}}
\newcommand{\wF}{\widetilde{F}}
\newcommand{\wf}{\widetilde{f}}
\newcommand{\wg}{\widetilde{g}}
\newcommand{\wh}{\widetilde{h}}
\newcommand{\wj}{\widetilde{j}}
\newcommand{\wb}{\widetilde{b}}
\newcommand{\wa}{\widetilde{a}}

\newcommand{\sC}{{\mathcal C}}

\newcommand{\sH}{{\mathcal H}}

\newcommand{\sO}{{\mathcal O}}

\newcommand{\sQ}{{\mathcal Q}}

\newcommand{\sU}{{\mathcal U}}

\newcommand{\sX}{{\mathcal X}}

\newcommand{\C}{{\mathbb C}}

\newcommand{\BP}{{\mathbb P}}

\def\Sym{\mathop{\rm Sym}\nolimits}

\def\Hom{\mathop{\rm Hom}\nolimits}

\title[Lagrangian loci in moduli of abelian surfaces]{Lagrangian loci in moduli of \\ abelian surfaces}
\author[Jun-Muk Hwang]{Jun-Muk Hwang} 
\address{Center for Complex Geometry,
Institute for Basic Science (IBS),
Daejeon 34126, Republic of Korea}
\email{jmhwang@ibs.re.kr}
\thanks{This work was supported by the Institute for Basic Science (IBS-R032-D1).}

\begin{document}

\begin{abstract} We show that any smooth surface germ in the moduli of abelian surfaces arises from a Lagrangian fibration of abelian surfaces. By  Donagi-Markman's cubic condition, the key issue of the proof is to find a suitable affine structure with a compatible cubic form on the base space of the family. We achieve this by analyzing the properties of cubic forms in two variables and proving the existence of the solution of the resulting partial differential equations by Cauchy-Kowalewski Theorem. Modifying the argument, we show also that a smooth curve germ in the moduli of abelian surfaces arises from a Lagrangian fibration if and only if the curve is a null curve with respect to the natural holomorphic conformal structure on the moduli of abelian surfaces.
\end{abstract}

\maketitle



\section{Introduction}
We study a problem on the local geometry of moduli of abelian varieties. As the problem is local, it is convenient to use Siegel upper half spaces instead of moduli of abelian varieties. Let us recall basic terminology.

Let $\sH_n$ be the Siegel upper half space of $n\times n$ complex symmetric matrices with positive definite imaginary part. Let $\Delta_n$ be the $n \times n$ integral diagonal matrix whose entries are the elementary divisors of a given polarization type (p.304 of \cite{GH}). Let $\varpi: \sX_n \to \sH_n$ be a universal abelian variety of the given polarization type, whose fiber over $Z \in \sH_n$ is $\C^n$ modulo the lattice generated by $2n$ columns of the $n \times 2n$ matrix $(\Delta_n, Z)$.
Let $\Sigma \subset \sH_n$ be the zero section of $\varpi$.

Our main interest is Lagrangian loci in $\sH_n$, i.e., submanifolds that can be realized as the images of the period maps of Lagrangian fibrations of abelian varieties. Let us define it more precisely.

\begin{definition}\label{d.Lagrangian}
A germ of complex submanifolds $s \in S \subset \sH_n$ is {\em a Lagrangian locus} if there exist a germ of complex manifolds $z \in B$ of dimension $n$ and a submersion $\phi: B \to S$ with $\phi (z) =s$ such that the fiber product $X := B \times_{\sH_n} \sX_n$  equipped with the induced morphism $\pi: X \to B$  and the induced zero section $\Sigma_X \subset X$ has
 a Lagrangian structure (in the sense of Section 1 of \cite{DM}), i.e.,  a holomorphic symplectic form $\omega$ on $X$ whose restrictions to  the submanifolds $\Sigma_X$ and $ X_y := \pi^{-1}(y)$ are zero for all $y \in B$.
\end{definition}

Note that the dimension of a Lagrangian locus $S$ is at most $n$. We raise the following question.

\begin{question}\label{q}
How do we characterize Lagrangian loci in $\sH_n$  in terms of the geometry of $\sH_n$? \end{question}

This question is closely related to the work \cite{DM}, where  Donagi and Markman studied which   families of abelian varieties admit   Lagrangian structures.
In Theorem 1 of \cite{DM}, they gave
 a necessary  condition, called `infinitesimal cubic condition', and in Theorem 2 of \cite{DM} (also see Theorem \ref{t.DM} below), they strengthened it to a necessary and sufficient condition, called `global cubic condition'.  The infinitesimal cubic condition is of algebraic nature and  can be translated directly to a necessary condition for
 Lagrangian locus in the setting of Definition \ref{d.Lagrangian}. However, the global cubic condition is formulated in terms of an affine structure on the base manifold of the family, which is a differential geometric condition. It is not immediate how to translate the latter into the language of the geometry of period loci in $\sH_n$. In this paper, we work out this translation when $n=2$ and give the following answer to Question \ref{q} when $n=2$.

 \begin{theorem}\label{t.main}
 Let $s \in S$ be a germ of complex submanifolds  in $\sH_2$. \begin{itemize}
 \item[(i)] If $\dim S = 2$, then $S$ is a Lagrangian locus.
     \item[(ii)] If $\dim S =1$, then $S$ is a Lagrangian locus if and only if it is a null curve with respect to the natural holomorphic conformal structure on $\sH_2$. \end{itemize} \end{theorem}

         The holomorphic conformal structure and  null curves in (ii) are explained in Definition \ref{d.conformal}.

 To prove Theorem \ref{t.main}, we translate Donagi-Markman's condition in terms of the geometry of $\sH_n$ and reduce the problem to the existence of a collection of closed forms on $B$ satisfying the cubic condition. To find such closed forms, we show that the  partial differential equations describing them are of Cauchy-Kowalewski type, hence have local analytic solutions.

 \medskip
  Theorem \ref{t.main} says that Donagi-Markman's {\em infinitesimal} cubic condition is sufficient to characterize Lagrangian loci in $\sH_2$, i.e., the integrability part of their global cubic condition can be always satisfied by a suitable choice of cubic forms.  It is natural to ask whether this still holds for $\sH_n, n>2$. Our methods can be adapted to certain types of germs of submanifolds in $\sH_n, n>2$, but there seem to be other types of germs which need new ideas. My guess is that the 2-dimensional case is exceptional and additional conditions appear in higher dimensions.

\section{Calculus of coframes in dimension 2}

\begin{definition}\label{d.coframe}
Fix  a complex vector space $V$.
Let $M$ be a complex manifold with $\dim M = \dim V$.
A {\em coframe } on $M$ is a $V$-valued 1-form $\theta$ on $M$ such that $\theta_y: T_y M \to V$ is an isomorphism for any $y \in M$. A coframe $\theta$ is {\em closed} if the $V$-valued 2-form ${\rm d} \theta$ is zero. Two coframes $\theta$ and $\wt$ are {\em linearly independent} at $y \in M$ if $\theta_y$ and $\wt_y$ are linearly independent  in $\Hom (T_y M, V)$. They are {\em pointwise linearly independent} if they are linearly independent at every point of $M$.  \end{definition}

\begin{proposition}\label{p.calculus}
Let $M$ be  a 2-dimensional complex manifold. Let
$\theta$ and $\wt$ be two coframes on $M$ which are linearly independent at a point $z \in M$. Then there exist a neighborhood $B \subset M$ of $z$ and two holomorphic functions $F$ and $\widetilde{F}$ on $B$ such that $F + \wF $ is nowhere vanishing and  $F \theta + \wF \wt$ is a closed coframe on $B$. \end{proposition}

\begin{proof}
Fix an identification $V = \C^2$ to write $$\theta =(\theta_1, \theta_2), \ \wt = (\wt_1, \wt_2)$$ for some holomorphic $1$-forms $\theta_1, \theta_2, \wt_1, \wt_2$ on $M$.
Let $(x, y)$ be a local holomorphic coordinate system in a neighborhood $O \subset M$ of $z$. We can find holomorphic functions $f, g, h, j, \wf, \wg, \wh, \wj$ on $O$ such that  $$
\left(\begin{array}{c} \theta_1 \\ \theta_2 \end{array} \right) =
 \left( \begin{array}{cc}
f & g \\ h & j \end{array} \right) \left( \begin{array}{c} {\rm d} x \\ {\rm d}y  \end{array} \right) \ \mbox{ and } \ \left(\begin{array}{c} \wt_1 \\ \wt_2 \end{array} \right) =
 \left( \begin{array}{cc}
\wf & \wg \\ \wh & \wj \end{array} \right) \left( \begin{array}{c} {\rm d} x \\ {\rm d} y  \end{array} \right). $$
 Since $\theta$ and $\wt$ are linearly independent at $z$, we may assume, by shrinking $O$ if necessary, that \begin{equation}\label{e.0}  f \wh - h \wf \mbox{ is nowhere zero on } O.                                         \end{equation}

We want to find holomorphic functions $F$ and $\wF$ on $O$ such that $\vartheta := F \theta + \wF \wt$ is a closed coframe on $O$. In terms of the components $\vartheta= (\vartheta_1, \vartheta_2),$ we need ${\rm d} \vartheta_1 = {\rm d} \vartheta_2 =0$. From $$ \vartheta_1 = (fF+\wf \wF)\ {\rm d}x + (g F + \wg \wF) \ {\rm d} y  $$ $$ \vartheta_2 = (h F + \wh \wF) \ {\rm d} x + (j F + \wj \wF) \ {\rm d} y,$$
the condition ${\rm d} \vartheta_1 =0$ reads
\begin{equation}\label{e.1} g \frac{\partial F}{\partial x} + \wg \frac{\partial \wF}{\partial x} - f \frac{\partial F}{\partial y} - \wf \frac{\partial \wF}{\partial y} + F(\frac{\partial g}{\partial x} - \frac{\partial f}{\partial y}) + \wF (\frac{\partial \wg}{\partial x} - \frac{\partial \wf}{\partial y}) = 0, \end{equation}
and the condition ${\rm d} \vartheta_2 = 0$ reads
\begin{equation}\label{e.2}
j \frac{\partial F}{\partial x} + \wj \frac{\partial \wF}{\partial x} - h \frac{\partial F}{\partial y} - \wh \frac{\partial \wF}{\partial y} + F(\frac{\partial j}{\partial x} - \frac{\partial h}{\partial y}) + \wF (\frac{\partial \wj}{\partial x} - \frac{\partial \wh}{\partial y}) = 0.\end{equation}
The two equations (\ref{e.1}) and (\ref{e.2}) can be combined into
$$
\left( \begin{array}{cc}
f & \wf \\ h & \wh \end{array} \right)
\left(\begin{array}{c}
\frac{\partial F}{\partial y} \\ \frac{\partial \wF}{\partial y} \end{array} \right) =
\left( \begin{array}{c} g \frac{\partial F}{\partial x} + \wg \frac{\partial \wF}{\partial x}  + F(\frac{\partial g}{\partial x} - \frac{\partial f}{\partial y}) + \wF (\frac{\partial \wg}{\partial x} - \frac{\partial \wf}{\partial y}) \\ j \frac{\partial F}{\partial x} + \wj \frac{\partial \wF}{\partial x}  + F(\frac{\partial j}{\partial x} - \frac{\partial h}{\partial y}) + \wF (\frac{\partial \wj}{\partial x} - \frac{\partial \wh}{\partial y}) \end{array} \right).
$$
By (\ref{e.0}), we can write it as
$$
\left(\begin{array}{c}
\frac{\partial F}{\partial y} \\ \frac{\partial \wF}{\partial y} \end{array} \right) =
 \left( \begin{array}{cc}
f & \wf \\ h & \wh \end{array} \right)^{-1} \left( \begin{array}{c} g \frac{\partial F}{\partial x} + \wg \frac{\partial \wF}{\partial x}  + F(\frac{\partial g}{\partial x} - \frac{\partial f}{\partial y}) + \wF (\frac{\partial \wg}{\partial x} - \frac{\partial \wf}{\partial y}) \\ j \frac{\partial F}{\partial x} + \wj \frac{\partial \wF}{\partial x}  + F(\frac{\partial j}{\partial x} - \frac{\partial h}{\partial y}) + \wF (\frac{\partial \wj}{\partial x} - \frac{\partial \wh}{\partial y}) \end{array} \right).
$$
This system  of partial differential equations is of Cauchy-Kowalewski type (e.g. Theorem 4 on p. 39 of \cite{Sp}). Thus for any holomorphic functions $F^o(x)$
and $\widetilde{F}^o(x)$ on the curve $(y=0) \subset O$,
we have holomorphic solutions  $F, \wF$ in a neighborhood $O' \subset O$ of $z$ such that $$F(x, 0) = F^o(x)\mbox{ and }\wF(x, 0) = \wF^o(x)$$ for all $(x,0) \in O'$.
Let us choose constant initial values $F^o \equiv 1$ and $\wF^o \equiv 0$.
Then for the solution $F$ and $\wF$  on $O'$ with these initial values, the $V$-valued closed 1-form $\vartheta = F \theta + \wF \wt$ remains a coframe in a neighborhood $B \subset O'$ of $z$ and $F+\wF$ is nowhere vanishing on $B$. This completes the proof.
\end{proof}

\section{Properties of cubic forms in 2 variables}

Here we collect some elementary properties of cubic forms in 2 variables.

\begin{definition}\label{d.chi}
Let $V$ be a $2$-dimensional complex vector space.
Define the injective homomorphism
$\alpha: \Sym^3 V^* \to \Hom (V, \Sym^2 V^*)$ by sending a cubic form $f(\cdot, \cdot, \cdot)$ to  $\alpha^f \in \Hom (V, \Sym^2 V^*)$ defined by $$v \in V \ \mapsto \ \alpha^f(v) \ = \  3 \ f(v, \cdot, \cdot) \ \in \Sym^2 V^*.$$ \begin{itemize} \item[(1)] A cubic form $f \in \Sym^3 V^*$ is said to be {\em degenerate} if
$$0 \neq {\rm Ker}(\alpha^f) = \{ u \in V, f(u, v, w) =0 \mbox{ for all } v, w \in V\}.$$ Otherwise, we say that it is {\em nondegenerate}.
\item[(2)] The set of degenerate cubic forms is exactly the affine cone $$\sC:= (\widehat{v_3(\BP V^*)} \setminus 0) \subset \Sym^3 V^*$$ of the third Veronese variety $v_3(\BP V^*) \subset \BP \Sym^3 V^*$. As ${\rm Im}(\alpha^f)$ has dimension 1 for any $f \in \sC,$ we have a natural morphism $$\kappa: \sC \to \BP \Sym^2 V^*$$ sending $[f]$ to the point $[{\rm Im}(\alpha^f)].$ It is clear that the image of $\kappa$ is the second Veronese variety $v_2(\BP V^*)$ and all fibers of $\kappa$ are biregular to $\C^*$.
\item[(3)] Let $\sU \subset (\Sym^3 V^* \setminus 0)$  be the set of nondegenerate cubic forms. As $\dim {\rm Im}(\alpha^f) =2$ for $f \in \sU$, we have a natural morphism $$\chi: \sU \to {\rm Gr}(\Sym^2 V^*; 2)$$
sending $f \in \sU$ to the point $[{\rm Im}(\alpha^f)]$ of the Grassmannian of $2$-dimensional subspaces in $\Sym^2 V^*$. \end{itemize}
\end{definition}

\begin{remark}
Because of the coefficient 3 in the definition of $\alpha^f$ in Definition \ref{d.chi}, we can compute $\alpha^f(v)$ by taking derivative of $f$ regarded as a polynomial function on $V$ by $v$ regarded as a vector field on $V$ (e.g. Exercise 3.6 in Chapter 1 of \cite{St}). We use this below in the computations of Propositions \ref{p.chi} and \ref{p.xi}. \end{remark}

\begin{proposition}\label{p.zeta}
In the setting of Definition \ref{d.chi} (2), let $W$ be a vector space of dimension $2$ and let $\varphi: W \to \Sym^2 V^*$ be a homomorphism with $\dim {\rm Im}(\varphi) =1 $ and $[{\rm Im}(\varphi)] \in v_2(\BP V^*) \subset \BP \Sym^2 V^*$.
For any cubic form $f \in \kappa^{-1}([{\rm Im}(\varphi)])$, we can find two linearly independent isomorphisms $\zeta: W \to V$ and $\widetilde{\zeta}: W \to V$ such that \begin{itemize} \item[(1)] $\alpha^f( \zeta(w))= \varphi(w) = \alpha^f(\widetilde{\zeta}(w))$ for all $ w \in W;$   and \item[(2)] for any $a, \widetilde{a} \in \C$ satisfying $a+\widetilde{a} \neq 0$, $\alpha^{\frac{1}{a+ \widetilde{a}} f} (a \zeta (w) + \widetilde{a} \widetilde{\zeta}(w)) = \varphi(w)$ for all $w \in W.$  \end{itemize}
 \end{proposition}

\begin{proof}
Choose a basis $\{ w_1, w_2\} $ of $W$ with $\varphi(w_1) =0$ and a basis $\{ v_1, v_2\}$ of $V$ with $\alpha^f(v_1) =0$. Then set $$\zeta(w_1) = v_1, \widetilde{\zeta}(w_1) = 2 v_1 \ \mbox{
and } \ \zeta(w_2) = \widetilde{\zeta}(w_2) = c v_2$$ where $c \in \C^*$ is determined by $ c \ \alpha^f(v_2)  = \varphi(w_2).$ Then $\zeta$ and $\widetilde{\zeta}$ satisfy all the required properties.  \end{proof}

\begin{proposition}\label{p.chi}
In Definition \ref{d.chi} (3), fix a basis $x, y $ of $V^*$ and identify $\Sym^2 V^*$ (resp. $\Sym^3 V^*$) with the space of homogeneous quadratic (resp. cubic) polynomials in $x, y$.
\begin{itemize}
\item[(1)] The natural ${\rm GL}(V)$-action on $\sU$ has  exactly two orbits in $\sU$,  represented by $x^3+ y^3$ and $x^2y,$ respectively.
    \item[(2)] The morphism $\hat{\chi}: \sU \to \wedge^2 (\Sym^2 V^*)$ defined by $$\hat{\chi}(f) = \frac{\partial f}{\partial x} \wedge \frac{\partial f}{\partial y} \in \wedge^2 (\Sym^2 V^*)$$
                       descends to $\chi: \sU \to {\rm Gr}(\Sym^2 V^*; 2)$  under the Pl\"ucker embedding  ${\rm Gr}(\Sym^2 V^*; 2) \subset \BP \wedge^2 (\Sym^2 V^*)$.
        \item[(3)] The ${\rm GL}(V)$-equivariant morphism $\chi$ sends $x^3+ y^3$ to $[x^2 \wedge y^2]$ with $$ \chi^{-1}([x^2 \wedge y^2]) = \{ c x^3 + e y^3, \ c, e \in \C^*\}$$ and sends $x^2y$ to $[xy \wedge x^2]$ with $$\chi^{-1}([xy \wedge x^2]) = \{ cx^3 + e x^2y, \ c \in \C, e \in \C^*\}.$$
            \item[(4)] The morphism $\chi: \sU \to {\rm Gr}(\Sym^2 V^*;2)$ is a submersion, i.e., is surjective and smooth.
            \end{itemize} \end{proposition}

\begin{proof}
We present only  the proof of the smoothness of the morphism $\chi$ in (4). The rest is straightforward.

To check that $\chi$ is a smooth morphism, it suffices to show that the differential of the morphism $\hat{\chi}: \sU \to \wedge^2(\Sym^2 V^*)$ in (2),  $${\rm d}_u \hat{\chi}: T_u \sU \to T_{\chi(u)} (\wedge^2 (\Sym^2 V^*)), $$ has rank 3  at each point $u \in \sU$.   By the ${\rm GL}(V)$-equivariance of $\chi$,
it is enough to check this when $u = x^3 + y^3$ and $u= x^2 y$.
Let us use the coordinates $(a, b, c, e)$ on $\sU \subset \Sym^3 V^*$ such that a cubic form $f \in \Sym^3 V^*$ is written as
$$f = a x^3 + b y^3 + c x^2y + e xy^2.$$
For $f \in U,$ the plane $\chi(f)$ of  $\wedge^2 V^*$ is spanned by $$\frac{\partial f}{\partial x} = 3a x^2 + 2 c xy + e y^2 \mbox{ and } \frac{\partial f}{\partial y} = 3 b y^2 + c x^2 + 2e xy.$$
Then the equation
$$ \frac{\partial f}{\partial x} \wedge \frac{\partial f}{\partial y} = (9 ab - ce) x^2 \wedge y^2 + (6ae -2c^2) x^2 \wedge xy + (2 e^2 -6bc) y^2 \wedge xy$$
says that in terms of the Pl\"ucker coordinates $(p, q, r)$ for ${\rm Gr}(\Sym^2 V^*; 2)$ given by the coefficients of $x^2 \wedge y^2, x^2 \wedge xy$ and $y^2 \wedge xy$, the morphism $\hat{\chi}: \sU \to \wedge^2(\Sym^2 V^*)$ is given by
$$ p= 9ab -ce, \ q = 6ae -2c^2, \ r = 2e^2 -6 bc.$$
The Jacobian matrix of $\hat{\chi}$
$$\left( \begin{array}{cccc}
9b & 9a & -e & -c \\
6e & 0 & -4c & 6a \\
0 & -6c & -6b & 4e \end{array} \right)$$
clearly has rank 3 at the points $u=x^3 + y^3$ and  $u = x^2y$. \end{proof}

\begin{proposition}\label{p.xi}
In the setting of Definition \ref{d.chi} (3) and Proposition \ref{p.chi}, let $W$ be a vector space of dimension $2$ and let $\varphi: W \to \Sym^2 V^*$ be an injective homomorphism with the corresponding point $[{\rm Im}(\varphi)]  \in {\rm Gr}(\Sym^2 V^*; 2)$.
\begin{itemize}
\item[(i)] Each element $f \in \chi^{-1} ([{\rm Im}(\varphi)])$  determines a unique isomorphism $\xi^f: W \to V$ satisfying $\alpha^f \circ \xi^f = \varphi.$
\item[(ii)] For two elements $f, \wf \in \chi^{-1}( [{\rm Im} (\varphi)]),$  there exists  $c \in \C^*$ satisfying $\xi^f = c \ \xi^{\wf}$ if and only if $c\ f = \wf.$
\item[(iii)] Let $f, \wf \in \chi^{-1}( [{\rm Im} (\varphi)])$ be  two linearly independent elements in $\Sym^3 V^*$ such that $a\ \xi^f + \wa \ \xi^{\wf} \in \Hom(W, V)$ is an isomorphism for some $a, \wa \in \C$. Then there are unique  $b, \wb \in \C$ such that $$b \ f + \wb \ \wf \in \chi^{-1}( [{\rm Im}(\varphi)]) \mbox{ and }
     \xi^{b f + \wb \wf} = a \ \xi^f + \wa \ \xi^{\wf}.$$ \end{itemize} \end{proposition}

\begin{proof}
When $f \in \chi^{-1} ([{\rm Im}(\varphi)])$, we have
${\rm Im}(\alpha^f) = {\rm Im} (\varphi)$. Thus there exists a unique isomorphism $\xi^f: W \to V$ satisfying (i).

 If $f, \wf \in \chi^{-1}( [{\rm Im} (\varphi)])$ satisfy $\xi^f = \xi^{\wf}$, then
$$\alpha^f \circ \xi^f = \varphi = \alpha^{\wf} \circ \xi^{\wf} = \alpha^{\wf} \circ \xi^f$$  Hence $\alpha^{(f - \wf)} \circ \xi^{f} =0,$ which implies $f - \wf \in {\rm Ker}(\alpha)$. Thus $f = \wf$ from the injectivity of $\alpha$.
It is clear that $\xi^{cf} = c^{-1} \ \xi^f$ for any $c \in \C^*$. Thus we have $\xi^f = c \ \xi^{\wf}$ for some $c \in \C^*$ if and only if $c\ f = \wf.$ This proves (ii).

As in the proof of Proposition \ref{p.chi}, by the ${\rm GL}(V)$-equivariance, we may check (iii) in the following  two special cases.

{\bf Case 1}. When $f = x^3+y^3.$  From Proposition \ref{p.chi} (3), we may put $\wf = c x^3 + e y^3 $ for some  $c \neq e \in \C^*$.
In terms of the dual basis $\frac{\partial}{\partial x}, \frac{\partial}{\partial y}$ of $V$, the isomorphism $\xi^{\wf}: W \to V$ satisfies \begin{equation}\label{e.xi1-1}
\frac{\partial}{\partial x} = \xi^{\wf} \circ \varphi^{-1} (\frac{\partial \wf}{\partial x}) = \xi^{\wf} \circ \varphi^{-1}(3c \ x^2) \end{equation} and \begin{equation}\label{e.xi1-11}
\frac{\partial}{\partial y}  = \xi^{\wf} \circ \varphi^{-1}(\frac{\partial \wf}{\partial y}) = \xi^{\wf} \circ \varphi^{-1}(3e \ y^2). \end{equation}
In the same way, we have \begin{equation}\label{e.xi1-2}
\frac{\partial}{\partial x} = \xi^{f} \circ \varphi^{-1}  (3 \ x^2) \end{equation} and \begin{equation}\label{e.xi1-21}
\frac{\partial}{\partial y} = \xi^{f} \circ \varphi^{-1}(3 \ y^2). \end{equation}
Putting $$f':= b \ f + \wb \ \wf = (b + \wb c) x^3 + (b + \wb e) y^3,$$ we obtain by a similar computation,
\begin{equation}\label{e.xi1-3}
\frac{\partial}{\partial x} = \xi^{f'} \circ \varphi^{-1}(3(b + \wb c) x^2) \end{equation} and \begin{equation}\label{e.xi1-31} \frac{\partial}{\partial y} = \xi^{f'} \circ \varphi^{-1}(3(b + \wb e) y^2). \end{equation}
We want to have \begin{equation*} (a \xi^f + \wa \xi^{\wf}) \circ \varphi^{-1}(x^2) = \xi^{f'} \circ \varphi^{-1}(x^2) \end{equation*} and \begin{equation*} (a \xi^f + \wa \xi^{\wf}) \circ \varphi^{-1}(xy) = \xi^{f'} \circ \varphi^{-1}(xy),\end{equation*} which can be translated by (\ref{e.xi1-1}) -- (\ref{e.xi1-31}) into the equations   \begin{equation*}
 \frac{1}{3(b + \wb c)} \frac{\partial}{\partial x} =  \frac{a}{3} \frac{\partial}{\partial x} + \frac{\wa}{3c} \frac{\partial}{\partial x} \end{equation*} and \begin{equation*}
\frac{1}{3(b + \wb e)} \frac{\partial}{\partial y} = \frac{a}{3} \frac{\partial}{\partial y} + \frac{\wa}{3e} \frac{\partial}{\partial y}. \end{equation*}
Thus it suffices to solve the two linear equations for $(b, \wb)$  $$b + c \wb = \frac{c}{ca + \wa}, \ \mbox{ and } \ b + e \wb = \frac{e}{ea + \wa}$$ where $ca + \wa \neq 0 \neq ea + \wa$ because of the assumption  that $a \xi^f + \wa \xi^{\wf}$ is an isomorphism. Since $c\neq e$, the system has a unique solution for $b, \wb$.

{\bf Case 2}. When $f = x^2y.$    From Proposition \ref{p.chi} (3), we may put $\wf = c x^3 + e x^2 y$ for some  $c, e \in \C^*$. We have
\begin{equation*}
\frac{\partial}{\partial x}  = \xi^{\wf} \circ \varphi^{-1}(\frac{\partial \wf}{\partial x}) = \xi^{\wf} \circ \varphi^{-1} (3c x^2 + 2e xy)\end{equation*} and \begin{equation*}  \frac{\partial}{\partial y} = \xi^{\wf} \circ  \varphi^{-1}(\frac{\partial \wf}{\partial y}) =   \xi^f \circ \varphi^{-1}(ex^2) \end{equation*}
 which give \begin{equation}\label{e.xi2-1}
\xi^{\wf} \circ \varphi^{-1}(x^2) = \frac{1}{e} \frac{\partial}{\partial y} \end{equation} and \begin{equation}\label{e.xi2-11}  \xi^{\wf} \circ \varphi^{-1}(xy) = \frac{1}{2e} \frac{\partial}{\partial x} - \frac{3c}{2e^2} \frac{\partial}{\partial y}.\end{equation}
Similarly, we have \begin{equation}\label{e.xi2-2}
\xi^{f} \circ \varphi^{-1}(x^2) =  \frac{\partial}{\partial y} \end{equation} and \begin{equation}\label{e.xi2-21} \xi^{f} \circ \varphi^{-1}(xy) = \frac{1}{2} \frac{\partial}{\partial x}.\end{equation}
Setting $f'= b f + \wb \wf = \wb c x^3 + (b + \wb e) x^2y$, we have by a similar computation, \begin{equation}\label{e.xi2-3}
\xi^{f'}\circ \varphi^{-1}(x^2) = \frac{1}{b + \wb e} \frac{\partial}{\partial y} \end{equation} and \begin{equation}\label{e.xi2-31} \xi^{f'} \circ \varphi^{-1}(xy) = \frac{1}{2(b+ \wb e)} \frac{\partial}{\partial x} - \frac{ 3\wb c}{2(b + \wb e)^2} \frac{\partial}{\partial y}. \end{equation}
We want to have \begin{equation*} (a \xi^f + \wa \xi^{\wf}) \circ \varphi^{-1}(x^2) = \xi^{f'} \circ \varphi^{-1}(x^2) \end{equation*} and \begin{equation*} (a \xi^f + \wa \xi^{\wf}) \circ \varphi^{-1}(xy) = \xi^{f'} \circ \varphi^{-1}(xy),\end{equation*} which can be translated by (\ref{e.xi2-1}) -- (\ref{e.xi2-31}) to the equations   \begin{equation*}a \frac{\partial}{\partial y} + \frac{\wa}{e} \frac{\partial}{\partial y} = \frac{1}{b + \wb e} \frac{\partial}{\partial y} \end{equation*} and \begin{equation*} \frac{a}{2} \frac{\partial}{\partial x} + \frac{\wa}{2e} \frac{\partial}{\partial x} - \frac{3 \wa c}{2e^2} \frac{\partial}{\partial y} = \frac{1}{2(b + \wb e)} \frac{\partial}{\partial x} - \frac{3 \wb c}{2(b + \wb e)^2} \frac{\partial}{\partial y}.\end{equation*} Since $ae + \wa \neq 0$ by the assumption that $a \xi^f + \wa \xi^{\wa}$ is an isomorphism and  $c \neq 0$, we obtain
$$ b + \wb e = \frac{e}{ae + \wa} \ \mbox{ and } \
e^2 \wb = \wa (b + \wb e)^2,$$
with the solutions
$$ b = \frac{a e^2}{(ae + \wa)^2} \ \mbox{ and } \  \wb = \frac{\wa}{(ae + \wa)^2}.$$
\end{proof}

\section{Proof of Theorem \ref{t.main}}

Let us recall `Global cubic condition' from Theorem 2 of \cite{DM}. To state it, we fix a vector space $V$ of dimension $n$ and regard $\sH_n$ as an open subset in $\Sym^2 V^*$ via a fixed identification $\C^n \cong V^*$. Note that  Donagi-Markman (in the paragraph before Theorem 2 of \cite{DM}) uses $\Sym^2 V$ in place of our $\Sym^2 V^*$. For us, it is more convenient to use the dual space.

\begin{theorem}\label{t.DM}
Let $z \in B$ be a germ of $n$-dimensional complex manifold
and $\phi: B \to S \subset \sH_n \subset \Sym^2 V^*$ be a submersion to a germ $s \in S$ as in Definition \ref{d.Lagrangian}. Then for $X := B \times_{\sH_n} \sX_n$,  the fibration $\pi: X \to B$ admits a Lagrangian structure, if and only if there exist
\begin{itemize}
\item[(1)] a $V$-valued closed coframe $\vartheta$ on $B;$ and
    \item[(2)] a $\Sym^3 V^*$-valued function $\psi$ on $B$ such that $$\alpha^{\psi(y)}( \vartheta(w) ) = {\rm d}_y \phi (w) \in T_{\phi(y)} \sH_n = \Sym^2 V^*$$ for all $y \in B$ and $w \in T_y B.$ \end{itemize}
        Here, the identification $T_{\phi (y)} \sH_n = \Sym^2 V^*$ comes from the inclusion $\sH_n \subset \Sym^2 V^*$ and $\alpha^{\psi(y)} \in \Hom(V, \Sym^2 V^*)$ is from Definition \ref{d.chi} applied to $f= \psi(y) \in
        \Sym^3 V^*$. \end{theorem}

        This is a direct translation of Theorem 2 of \cite{DM} in our terminology from previous sections.
         In particular, affine structures used in \cite{DM} can be derived from closed coframes as follows.
If $\vartheta$ is a closed coframe on the complex manifold $B$, then Poincare lemma gives the affine holomorphic coordinate system $(x^1, \ldots, x^n),$ after replacing $B$ by a neighborhood of  $z \in B$, such that
$$\vartheta = {\rm d}x^1 \cdot v_1 + \cdots + {\rm d} x^n \cdot v_n$$ for some basis $\{ v_1, \ldots, v_n\}$ of $V$. These coordinates determine an affine structure on $B$.

\medskip
We use Theorem \ref{t.DM} for $n=2$ to prove Theorem \ref{t.main}.

\begin{proof}[Proof of Theorem \ref{t.main}(i)]
When $S \subset \sH_2$ is a surface germ, choose an embedding $\phi: B \to \sH_2 \subset \Sym^2 V^*$ with image $S$.
Define the holomorphic map $\phi_{*}: B \to {\rm Gr}(\Sym^2 V^*; 2)$ by $$ \phi_{*}(y) := [{\rm d}_y \phi(T_y B)] \mbox{ for each } y \in B,$$ where ${\rm d}_y \phi: T_y B \to T_{\phi(y)} \sH_2 = \Sym^2 V^*$ is the differential of $\phi$ at $y$.
By Proposition \ref{p.chi} (4), we can find two holomorphic maps $q, \wq : B \to \sU$ such that $\chi \circ q = \phi_* = \chi \circ \wq $ and  $q(y)$ and $\wq (y)$ are linearly independent in $\Sym^3 V^*$ for each $y \in B.$

By
Proposition \ref{p.xi} (i) applied to $W = T_y B$ for each $y \in B$,
with $\varphi = {\rm d}_y \phi : W=T_y B \to \Sym^2 V^*,$ we have  two $V$-valued coframes $\theta, \wt$ on $B$ defined by
$$\theta (w) := \xi^{q(y)} (w) \mbox{ and } \wt(w) :=
\xi^{\widetilde{q}(y)} (w) $$ for $w \in W= T_y B$ for each $y \in B$.  By Proposition \ref{p.xi} (ii), the two coframes $\theta$ and $\wt$ are pointwise linearly independent. Thus Proposition \ref{p.calculus} gives holomorphic functions $F, \wF$ on $B$ such that $\vartheta := F \theta + \wF \wt$ is a closed coframe on $B$.

To prove Theorem \ref{t.main} (i), it suffices to show, by Theorem \ref{t.DM}, that there are holomorphic functions $h, \wh$ on $B$ such that the $\Sym^3 V^*$-valued function on $B$ defined by $\psi = h \ q + \wh \widetilde{q}$ has image in $\sU$ and  satisfies
$$ \alpha^{\psi(y)}(\vartheta(w)) = {\rm d} \phi (w) $$ for all $y \in B$ and $w \in T_y B.$ In other words, we need to find $h, \wh$ satisfying
\begin{equation*} \xi^{h(y) q(y) + \wh(y) \wq (y)} = F(y) \xi^{q(y)} + \wF(y) \xi^{\wq(y)} \end{equation*}
for all $y \in B$.
This is equivalent to $$
\alpha^{h q + \wh \wq} (F \xi^{q}(w) + \wF \xi^{\wq}(w))  = {\rm d} \phi (w)$$  for all $w \in TB$.
Fixing two pointwise linearly independent sections $w_1$ and $w_2$ of $TB$ by shrinking $B$ if necessary, we can rewrite the above equation as two equations  in $\Sym^2 V^*$
\begin{equation}\label{e.th}
 \alpha^q(F \xi^q (w_i) + \wF \xi^{\wq} (w_i)) h +  \alpha^{\wq}(F \xi^q (w_i) + \wF \xi^{\wq} (w_i)) \wh =   {\rm d} \phi (w_i)
\end{equation}
for $i=1,2$. After fixing a basis of $\Sym^2 V^*$, the equations (\ref{e.th}) becomes
 a system of six inhomogeneous linear equations in two unknown functions $h$ and $\wh$ with coefficients in the ring $\sO (B)$ of holomorphic functions on $B$. By Proposition \ref{p.xi} (iii), this system of equations evaluated at each $y \in B$ has a unique solution. Thus (\ref{e.th}) has a unique solution in $\sO(B)$.
\end{proof}

Now we give a precise definition of the terms in Theorem \ref{t.main} (ii).

\begin{definition}\label{d.conformal}
When $\dim V=2$, the manifold $M = \Sym^2 V^*$  has a natural holomorphic conformal structure, i.e., a nondegenerate quadratic form determined up to  scalars on each tangent space $T_x M, x \in M$.
In fact,  the conic curve $v_2(\BP V^*) \subset \BP \Sym^2 V^*$ determines a conic bundle $\sQ\subset \BP TM$  and the quadratic form on $T_x M$ is just the equation of $\sQ_x \subset \BP T_x M$, determined up to scalars. A germ $S$ of curves in $\sH_2$ is called a {\em null curve} if $\BP T_x S \in \sQ_x$ for each $x \in S$. \end{definition}

\begin{proof}[Proof of Theorem \ref{t.main}(ii)]
If a curve germ $S \subset \sH_2$ is a Lagrangian locus,  the only-if part of Theorem \ref{t.DM} says that for each $x \in S$,
$$[T_x S] \in \BP T_x \sH_2 = \BP \Sym^2 V^*$$ is in the image of the morphism $\kappa$ in Definition \ref{d.chi} (2). Thus  $S$ is a null curve.

Let us prove the converse. The argument is parallel to the proof of Theorem \ref{t.main} (i). For a  germ $s \in S$ of null curves, choose a submersion $$\phi: B \to S \subset \sH_2 \subset \Sym^2 V^*$$ sending a germ of $2$-dimensional manifolds $z \in B$ to $s \in S$.
Define the holomorphic map $\phi_{*}: B \to v_2(\BP V^*) $ by $$ \phi_{*}(y) := [{\rm d}_y \phi (T_y B)] \mbox{ for each } y \in B,$$ where ${\rm d}_y \phi: T_y B \to T_{\phi(y)} \sH_2 = \Sym^2 V^*$ is the differential of $\phi$ at $y$.
By Definition \ref{d.chi} (ii), we can find a holomorphic map $q : B \to \sC$ such that $\kappa \circ q = \phi_*$.

By
Proposition \ref{p.zeta} applied to $W = T_y B$ for each $y \in B$ and $\varphi = {\rm d}_y \phi : W=T_y B \to \Sym^2 V^*,$ we have  two pointwise independent $V$-valued coframes $\theta, \wt$ on $B$ such that
$$ \alpha^{q(y)} (\theta (w)) = {\rm d}_y \phi (w) = \alpha^{q(y)}(\wt (w))$$ for all $y \in B$ and $w \in  T_y B$.  By Proposition \ref{p.calculus}, we have holomorphic functions $F, \wF$ on $B$ such that $F+ \wF \neq 0$ and  $\vartheta := F \theta + \wF \wt$ is a closed coframe on $B$ , after shrinking $B$ if necessary.

To apply Theorem \ref{t.DM}, we need a $\Sym^3 V^*$-valued holomorphic function $\psi$ on $B$ satisfying
$$ \alpha^{\psi (y)}(\vartheta(w)) = {\rm d}_y \phi (w) $$ for all $y \in B$ and $w \in T_y B.$
By Proposition \ref{p.zeta}, we can just choose $f= \frac{1}{F+ \wF} q$. \end{proof}


\end{document}